\documentclass{amsart}
\usepackage{amssymb,amscd,amsthm}
\usepackage{latexsym}
\usepackage{graphicx}
\date{Summer Semester 2008}

\newcommand{\Z}{{\mathbb Z}}
\newcommand{\R}{{\mathbb R}}
\newcommand{\C}{{\mathbb C}}

\newcommand{\T}{{\mathbb T}}

\newtheorem{theorem}{Theorem}

\newtheorem{lemma}[theorem]{Lemma}

\newtheorem{definition}[theorem]{Definition}
\begin{document}

\title[The Spectrum of the Almost Mathieu Operator]{{\rm Lecture
    series in the CRC 701\\ Summer Semester 2008}\\[5mm] {\Large The Spectrum of the Almost Mathieu Operator}}

\author{David Damanik}

\address{Department of Mathematics, Rice University, Houston, TX~77005, USA}

\email{damanik@rice.edu}

\maketitle

\begin{abstract}
These notes are based on a series of six lectures, given during my
stay at the CRC~701 in June/July 2008. The lecture series intended
to give a survey of some of the results for the almost Mathieu
operator that have been obtained since the early 1980's.
Specifically, the metal-insulator transition is discussed in
detail, along with its relation to the ten Martini problem via
duality and reducibility.
\end{abstract}

\section{Introduction and Overview}

\subsection{The Operator and the Main Results}

We will study the almost Mathieu operator
$$
[H^{\lambda,\alpha}_\omega \psi](n) = \psi(n+1) + \psi(n-1) + 2
\lambda \cos (2 \pi (\omega + n \alpha)) \psi(n).
$$

The potential $2 \lambda \cos (2 \pi (\omega + n
\alpha))$ is periodic if $\lambda = 0$ or $\alpha$ is rational and hence we will
only consider the case where $\lambda \not= 0$ and $\alpha$ is
irrational. Furthermore, by periodicity of the cosine, we consider
$\alpha$ and $\omega$ as elements of $\T = \R / \Z$. Finally, it is also easy
to see that $H^{\lambda,\alpha}_\omega =
H^{-\lambda,\alpha}_{\omega + \frac{1}{2}}$. Thus, we will only consider:
$\lambda > 0$, $\alpha \in \T$ irrational, $\omega \in
\T$.

The two theorems below concern the metal-insulator transition and
the ten Martini problem, that is, an exact description of the
spectral type of the operator, exhibiting a remarkably strict
transition at $\lambda = 1$, and the Cantor structure of the
spectrum. They are stated in the generality in which they are
currently known and summarize the results of many authors,
obtained over the course of about three decades.

\begin{theorem}[Metal-Insulator Transition]
{\rm (a)} If $\lambda < 1$, then for every $\alpha$ and every $\omega$, the spectrum is purely absolutely continuous.\\
{\rm (b)} If $\lambda = 1$, then for every $\alpha$ and all but countably many $\omega$, the spectrum is purely singular continuous.\\
{\rm (c)} If $\lambda > 1$, then for almost every $\alpha$ and almost every $\omega$, the spectrum is pure point and the eigenfunctions decay exponentially. \\
{\rm (d)} If $\lambda > 1$, then for generic $\alpha$ and every $\omega$, the spectrum is purely singular continuous.\\
{\rm (e)} If $\lambda > 1$, then for every $\alpha$ and generic $\omega$, the spectrum is purely singular continuous.
\end{theorem}

\begin{theorem}[Ten Martini Problem]
The spectrum of $H^{\lambda,\alpha}_\omega$ is a Cantor set, that is, it is closed and it contains no isolated points and no intervals.
\end{theorem}

We will present many of the main ideas that go into the proof of
these theorems. Some of the statements above will be proved
completely here, while the proof of others will only be sketched.
In the next subsection, we state the results we discuss in more
depth in subsequent sections.

\subsection{A Quick Guide to Proving the Main Results}

Consider the Hilbert space $L^2(\T \times
\Z)$ and the operator $H^{\lambda,\alpha} : L^2(\T \times \Z) \to L^2(\T \times \Z)$
given by
$$
[H^{\lambda,\alpha} \varphi](\omega,n) = \varphi(\omega,n+1) +
\varphi(\omega,n-1) + 2 \lambda \cos (2 \pi (\omega + n \alpha))
\varphi(\omega,n).
$$
Introduce the duality transform
$\mathcal{A} : L^2(\T \times \Z) \to L^2(\T \times \Z)$,
which is given by
$$
[\mathcal{A} \varphi](\omega,n) = \sum_{m \in \Z} \int_\T e^{-2\pi
i (\omega + n \alpha)m} e^{-2\pi i n \eta} \varphi(\eta,m) \,
d\eta.
$$
This definition assumes initially that $\varphi$ is such that the
sum in $m$ converges, but note that in terms of the Fourier
transform on $L^2(\T \times \Z)$, we have
$[\mathcal{A} \varphi](\omega,n) = \hat \varphi (n, \omega + n
\alpha)$, which may be used to extend the definition to all of $L^2(\T
\times \Z)$ and shows that $\mathcal{A}$ is unitary.

\begin{theorem}[Gordon-Jitomirskaya-Last-Simon
1997]
Suppose $\lambda > 0$ and $\alpha \in \T$ is irrational.\\
{\rm (a)} We have $H^{\lambda,\alpha} \mathcal{A} = \lambda \mathcal{A}
H^{\lambda^{-1},\alpha}$.
\\
{\rm (b)} If $H^{\lambda,\alpha}_\omega$ has
pure point spectrum for almost every $\omega \in \T$, then
$H^{\lambda^{-1},\alpha}_\omega$ has purely absolutely continuous
spectrum for almost every $\omega \in \T$.
\\
{\rm (c)} If $H^{\lambda,\alpha}_\omega$ has
some point spectrum for almost every $\omega \in \T$, then
$H^{\lambda^{-1},\alpha}_\omega$ has some absolutely continuous
spectrum for almost every $\omega \in \T$.
\end{theorem}

\begin{definition}
Fix $\lambda$ and $\alpha$. For $E \in \C$, the Lyapunov exponent is given by
$$
\gamma(E) = \inf_{n \ge 1} \frac1n \, \int \log \| M_E(n,\omega) \|
\, d\omega = \lim_{n \to \infty} \frac{1}{n} \log \| M_E(n,\omega) \| \;
\text{ for } \mu-\text{almost every } \omega \in \Omega,
$$
where
$$
M_E(n,\omega) = T_E(n,\omega) \times \cdots \times T_E(1,\omega),
\quad T_E(m,\omega) =
\begin{pmatrix} E - 2\lambda \cos (2 \pi (\omega + m \alpha)) & -1 \\ 1 & 0 \end{pmatrix}.
$$
\end{definition}

\begin{theorem}[Herman 1983]
The Lyapunov exponent obeys $\gamma(E) \ge \log \lambda$.
\end{theorem}

\begin{definition}\label{d.liou}
An irrational number $\alpha \in \T$ is called Liouville if there is a sequence of rational numbers $\frac{p_k}{q_k}$ with $q_k \to \infty$ such that
$$
\left| \alpha - \frac{p_k}{q_k} \right| < k^{-q_k}.
$$
\end{definition}

\begin{theorem}[Avron-Simon 1982]
Suppose $\alpha \in \T$ is Liouville. Then, for every $\lambda$ and $\omega$,
$H_\omega^{\lambda,\alpha}$ has purely continuous spectrum.
\end{theorem}

\begin{theorem}[Jitomirskaya-Simon 1994]
For every $\lambda$ and $\alpha$,
$H_\omega^{\lambda,\alpha}$ has purely continuous spectrum for generic $\omega$.
\end{theorem}

\begin{definition}
An irrational number $\alpha \in \T$ is called Diophantine if
there are constants $c = c(\alpha) > 0$ and $r = r(\alpha) > 1$
such that
$$
| \sin ( 2 \pi n \alpha ) | > \frac{c}{|n|^r} \quad \text{ for
every } n \in \Z \setminus \{0\}.
$$
Given such an $\alpha$, $\omega \in \T$ is called resonant if the
relation
$$
\left| \sin \left( 2\pi \left( \omega + \frac{n}{2} \alpha \right)
\right) \right| < \exp \left( -|n|^\frac{1}{2r} \right)
$$
holds for infinitely many $n \in \Z$; otherwise $\omega$ is called
non-resonant.
\end{definition}

It is known that Lebesgue almost every $\alpha$ is Diophantine.
Moreover, the set of resonant $\omega$'s is a dense $G_\delta$ set
(as can be seen directly from the definition) of zero Lebesgue
measure (by Borel-Cantelli).

\begin{theorem}[Jitomirskaya 1999]\label{t.jitoloc}
Suppose $\lambda > 1$, $\alpha \in \T$ is Diophantine, and $\omega
\in \T$ is non-resonant. Then, the almost Mathieu operator
$H_\omega^{\lambda,\alpha}$ has pure point spectrum with
exponentially decaying eigenfunctions.
\end{theorem}

\begin{theorem}[Puig~2004]\label{t.puigcantor}
Suppose $\alpha$ is Diophantine and $\lambda \not= 1$. Then,
$\Sigma^{\lambda,\alpha}$ is a
Cantor set.
\end{theorem}


\section{The Herman Estimate and Aubry Duality}

\subsection{The Herman Estimate}

The spectral type of $H_\omega^{\lambda,\alpha}$ can be studied by
looking at the solutions of the time-independent Schr\"odinger
equation:
$$
u(n+1) + u(n-1) + 2 \lambda \cos (2 \pi (\omega + n \alpha)) u(n) = E u(n).
$$
Notice that $u$ solves this equation if and only if
$$
\begin{pmatrix} u(n+1) \\ u(n) \end{pmatrix} = M_E(n,\omega)
\begin{pmatrix} u(1) \\ u(0) \end{pmatrix},
$$
where the transfer matrix $M_E(n,\omega)$ is given (at least for $n \ge 1$) by
$$
M_E(n,\omega) = T_E(n,\omega) \times \cdots \times T_E(1,\omega),
\quad T_E(m,\omega) =
\begin{pmatrix} E - 2\lambda \cos (2 \pi (\omega + m \alpha)) & -1 \\
  1 & 0 \end{pmatrix}.
$$
Here we leave the dependence on $\lambda$ and $\alpha$ implicit as
only $\omega$ will be varied in this section. Thus, decay or growth of
solutions is closely related to growth (of the norm) of the transfer
matrices. To measure the growth on an exponential scale, one
introduces the Lyapunov exponent. Initially, consider an
$\omega$-average and define
$$
\gamma(E) = \lim_{n \to \infty} \frac1n \, \int_\T \log \|
M_E(n,\omega) \| \, d\omega.
$$
Clearly, $\gamma(E) \ge 0$ since the matrices $M_E(n,\omega)$ have
determinant one and hence norm at least one. The existence of the
limit defining $\gamma(E)$ follows from Kingman's subadditive
ergodic theorem. In fact, this theorem also shows that
$$
\gamma(E) = \lim_{n \to \infty} \frac1n \, \log \| M_E(n,\omega)
\|
$$
for almost every $\omega \in \T$. For every such $\omega$, Osceledec'
theorem then shows that if $\gamma(E) > 0$, there is a one-dimensional
subspace of $(u(1),u(0))^T$ for which the norm of $(u(n+1),u(n))^T$
decays like $e^{-\gamma(E)n}$ as $n \to \infty$ and all linearly
independent initial conditions yield $e^{\gamma(E)n}$ growth.

\begin{theorem}[Herman 1983]
We have $\gamma(E) \ge \log \lambda$ for every $E$.
\end{theorem}

\begin{proof}
Setting $w = e^{2 \pi i\omega}$, we see that
$$
2 \lambda
\cos(2\pi (\omega + m \alpha)) = \lambda \left( e^{2 \pi i \alpha m} w
  + e^{-2 \pi i
\alpha m} w^{-1} \right) .
$$
Thus, the one-step transfer matrices have the form
$$
T_E(m, \omega) = \left( \begin{array}{cr} E - \lambda \left( e^{2
\pi i \alpha m} w + e^{-2 \pi i \alpha m} w^{-1} \right) & -1
\\ 1 & 0 \end{array} \right)
$$
If we define
$$
N_n(w) = w^n M_E(n,\omega) = (w T_E(n, \omega)) \cdots (w
T_E(1, \omega)),
$$
initially on $|w| = 1$, we see that $N_n$ extends to an entire
function and hence $w \mapsto \log \| N_n(w) \|$ is subharmonic.
Thus,
$$
\int_0^1 \log \| N_n(e^{2 \pi i\omega}) \| \, d\omega \ge \log \|
N_n(0) \| = n \log \lambda.
$$
Moreover, $\|N_n(e^{2 \pi i\omega})\| = \|M_E(n,\omega)\|$. Thus,
\begin{align*}
\gamma(E) & = \lim_{n \to \infty} \frac1n \int_\T \log \|
M_E(n,\omega) \| \, d\omega \\
& = \lim_{n \to \infty} \frac1n \int_\T \log \| N_n(e^{2 \pi i\omega})
\| \, d\omega \\
& \ge \log \lambda ,
\end{align*}
as claimed.
\end{proof}

\subsection{Aubry Duality}

The Herman estimate suggests strongly that there are exponentially
decaying solutions of the time-independent Schr\"odinger equation when
$\lambda > 1$. Notice, however, that we can infer this only for almost
every $\omega$ for any given $E$ and, moreover, we have treated only
exponential decay near $+\infty$. One may perform a similar analysis
near $-\infty$, but even if one finds two initial conditions which
yield exponential decays at $\pm \infty$, it is not clear whether they
coincide and give rise to a solution $u$ which decays at both ends and
hence is an exponentially decaying eigenfunction corresponding to the
eigenvalue $E$.

Nevertheless, it is indeed sometimes possible to show the
existence of genuine eigenfunctions $u$, this will be discussed in
Section~4. Aubry duality then constructs from such solutions
corresponding solutions for the dual coupling constant
$\lambda^{-1}$ and the dual energy $\lambda^{-1} E$. The special
form of these dual solutions will suggest that the spectral
measures for $\lambda^{-1}$ should have some absolutely continuous
component.

Let us demonstrate how this works. The relevant equations are
\begin{equation}\label{f.amoeveg}
u(n+1) + u(n-1) + 2 \lambda \cos(2 \pi (\omega + n \alpha)) u(n) =
E u(n),
\end{equation}
and the dual difference equation
\begin{equation}\label{f.damoeveg}
\tilde u(n+1) + \tilde u(n-1) + 2 \lambda^{-1} \cos(2 \pi (\tilde
\omega + n \alpha)) \tilde u(n) = (\lambda^{-1} E) \tilde u(n).
\end{equation}

\begin{lemma}\label{l.poinwisedual}
{\rm (a)} Suppose $u \in \ell^1(\Z)$ is a solution of
\eqref{f.amoeveg}. Consider its Fourier transform
$$
\hat u(\theta) = \sum_{m \in \Z} u(m) e^{2 \pi i m \theta}.
$$
Then, given any $\tilde \omega \in \T$, the sequence $\tilde u$
defined by
\begin{equation}\label{f.dualsoldef}
\tilde u(n) = \hat u(\tilde \omega + n \alpha) e^{2 \pi i n
\omega}
\end{equation}
is a solution of \eqref{f.damoeveg}.

{\rm (b)} Suppose $u \in \ell^2(\Z)$ is a solution of
\eqref{f.amoeveg}. Then, for $\tilde \omega$ from a full-measure
subset of $\T$, the sequence $\tilde u$ defined by
\eqref{f.dualsoldef} is a solution of \eqref{f.damoeveg}.
\end{lemma}

\begin{proof}
(a) If $u \in \ell^1(\Z)$, $\hat u \in C(\T)$ and we can evaluate
it pointwise. We have
\begin{align*}
&(\lambda^{-1} E) \tilde u(n) \\& = (\lambda^{-1} E) \hat u(\tilde
\omega + n \alpha) e^{2 \pi i n
\omega} \\
& = (\lambda^{-1} E) \sum_{m \in \Z} u(m) e^{2 \pi i m (\tilde
\omega + n \alpha)}
e^{2 \pi i n \omega}\\
& = \lambda^{-1} \sum_{m \in \Z} (E u(m)) e^{2 \pi i m (\tilde
\omega + n \alpha)}
e^{2 \pi i n \omega} \\
& = \lambda^{-1} \sum_{m \in \Z} \left[ u(m+1) + u(m-1) + 2
\lambda \cos (2 \pi (\omega + m \alpha)) u(m) \right] e^{2 \pi i m
(\tilde \omega + n \alpha)}
e^{2 \pi i n \omega} \\
& = \lambda^{-1} \sum_{m \in \Z} \left[ u(m+1) + u(m-1) + \lambda
\left( e^{2 \pi i (\omega + m \alpha)} + e^{-2 \pi i (\omega + m
\alpha)} \right) u(m) \right] e^{2 \pi i m \tilde \omega} e^{2 \pi i n
(\omega + m \alpha)}\\
& = \lambda^{-1} \sum_{m \in \Z} u(m+1) e^{2 \pi i m \tilde
\omega} e^{2 \pi i n (\omega + m \alpha)} + \lambda^{-1} \sum_{m
\in \Z} u(m-1) e^{2 \pi i m \tilde \omega} e^{2 \pi i n (\omega + m \alpha)} \\
& \qquad + \sum_{m \in \Z} u(m) e^{2 \pi i m \tilde \omega} e^{2
\pi i (n+1) (\omega + m \alpha)} +
\sum_{m \in \Z} u(m) e^{2 \pi i m \tilde \omega} e^{2 \pi i (n-1)
  (\omega + m \alpha)}  \\
& = \lambda^{-1} \sum_{m \in \Z} u(m) e^{2 \pi i (m-1) \tilde
\omega} e^{2 \pi i n (\omega + (m-1) \alpha)} + \lambda^{-1}
\sum_{m \in \Z} u(m) e^{2 \pi i (m+1) \tilde \omega} e^{2 \pi i n
  (\omega + (m+1) \alpha)} \\
& \qquad + \sum_{m \in \Z} u(m) e^{2 \pi i m \tilde \omega} e^{2
\pi i (n+1) (\omega + m \alpha)} +
\sum_{m \in \Z} u(m) e^{2 \pi i m \tilde \omega} e^{2 \pi i (n-1)
  (\omega + m \alpha)}  \\
& = \lambda^{-1} \left( e^{- 2 \pi i (\tilde \omega + n \alpha)} +
e^{2 \pi i (\tilde \omega + n \alpha)} \right) \sum_{m \in \Z}
u(m) e^{2 \pi i m (\tilde \omega + n \alpha)}
e^{2 \pi i n \omega} \\
& \qquad + \sum_{m \in \Z} u(m) e^{2 \pi i m (\tilde \omega +
(n+1) \alpha)} e^{2 \pi i (n+1) \omega} + \sum_{m \in \Z} u(m)
e^{2 \pi i m (\tilde \omega + (n-1) \alpha)} e^{2 \pi i (n-1) \omega}  \\
& = 2 \lambda^{-1} \cos(2 \pi (\tilde \omega + n \alpha)) \tilde
u(n) + \tilde u(n+1) + \tilde u(n-1).
\end{align*}

(b)  If $u \in \ell^2(\Z)$, $\hat u$ exists as an element of
$L^2(\T)$ and hence it is determined almost everywhere. Consider
$\tilde \omega$ from the full measure set of elements for which
all the quantities in the calculation above are determined. Then
carry out the calculation to verify that $\tilde u$ is indeed a
solution of \eqref{f.damoeveg} for the $\tilde \omega$ in
question.
\end{proof}

By pursuing these relations further, it is possible to show the following:

\begin{theorem}[Gordon-Jitomirskaya-Last-Simon
1997]\label{t.gjlsdual}
{\rm (a)} If $H^{\lambda,\alpha}_\omega$ has
pure point spectrum for almost every $\omega \in \T$, then
$H^{\lambda^{-1},\alpha}_\omega$ has purely absolutely continuous
spectrum for almost every $\omega \in \T$.
\\
{\rm (b)} If $H^{\lambda,\alpha}_\omega$ has
some point spectrum for almost every $\omega \in \T$, then
$H^{\lambda^{-1},\alpha}_\omega$ has some absolutely continuous
spectrum for almost every $\omega \in \T$.
\end{theorem}

Another consequence of Aubry duality is a formula relating the spectra of $H^{\lambda,\alpha}_\omega$ and
$H^{\lambda^{-1},\alpha}_\omega$. Note first that for $\alpha$ irrational, the spectrum of $H^{\lambda,\alpha}_\omega$ is independent of $\omega$ and may therefore be denoted by $\Sigma^{\lambda,\alpha}$. This follows from minimality of irrational rotations and strong operator convergence.

\begin{theorem}[Avron-Simon 1983]
We have $\Sigma^{\lambda,\alpha} = \lambda \Sigma^{\lambda^{-1},\alpha}$. In particular, $\Sigma^{\lambda,\alpha}$ is a Cantor set if and only if $\Sigma^{\lambda^{-1},\alpha}$ is a Cantor set.
\end{theorem}


\section{Exceptional Frequencies and Phases}

\subsection{Liouville Frequencies and the Gordon Method}

For the almost Mathieu operator $H_\omega^{\lambda,\alpha}$ with
$\lambda > 0$, $\alpha \in \T$ irrational, and $\omega \in \T$, we
investigate the spectral type by looking at the solutions of the
time-independent Schr\"odinger equation:
$$
u(n+1) + u(n-1) + 2 \lambda \cos (2 \pi (\omega + n \alpha)) u(n) = E u(n).
$$
We have seen in Section~2 that for $\lambda > 1$, the Lyapunov
exponent $\gamma(E)$ is uniformly bounded away from zero. By
Osceledec' theorem, it follows that for every $E$, there is an
$E$-dependent full measure set of $\omega$'s for which the
equation has an exponentially decaying solution at $+\infty$.
Similarly, there is also an exponentially decaying solution at
$-\infty$. Since it is known that almost everywhere with respect
to any spectral measure, there are polynomially bounded solutions,
this suggests that we should expect that there are exponentially
decaying eigenfunctions for spectrally almost every energy and
hence localization. However, here we need to interchange the
quantifiers and the application of Fubini then loses track of sets
of Lebesgue measure zero. Thus, the argument is inconclusive and
in fact wrong in general as we will see in this section.

The expected localization result fails when $\alpha$ is very well
approximated by rational numbers.

\begin{theorem}[Avron-Simon 1982]\label{t.avsim}
Suppose $\alpha \in \T$ is Liouville. Then, for every $\lambda$ and $\omega$,
$H_\omega^{\lambda,\alpha}$ has purely continuous spectrum.
\end{theorem}

The heart of the argument is the Cayley-Hamilton Theorem, which
for $\mathrm{SL}(2,\C)$ matrices $M$ takes the form
\begin{equation}\label{f.cayley-hamilton}
M^2 - \mathrm{Tr} M \cdot M + I = 0.
\end{equation}
Recall that the transfer matrices belong to $\mathrm{SL}(2,\C)$.
We will apply \eqref{f.cayley-hamilton} to these matrices when
there are suitable local repetitions. Explicitly, the following
version of Gordon's Lemma implements this.

\begin{lemma}\label{l.gordon}
Suppose $V : \Z \to \R$ obeys $V(n+p) = V(n)$ for some $p \in
\Z_+$ and $-p + 1 \le n \le p$, $E \in \R$, and $u$ solves
$$
u(n+1) + u(n-1) + V(n) u(n) = E u(n).
$$
Then, we have
\begin{equation}\label{f.gordonest}
\max \left\{ \left\| \begin{pmatrix} u(-p+1) \\ u(-p)
\end{pmatrix} \right\| , \left\| \begin{pmatrix} u(p+1) \\ u(p)
\end{pmatrix} \right\| , \left\| \begin{pmatrix} u(2p+1) \\ u(2p)
\end{pmatrix} \right\| \right\} \ge \frac12 \left\| \begin{pmatrix} u(1) \\ u(0)
\end{pmatrix} \right\|.
\end{equation}
\end{lemma}

\begin{proof}
By assumption, we have
$$
\begin{pmatrix} u(2p+1) \\ u(2p)
\end{pmatrix} = M_E(2p) \begin{pmatrix} u(1) \\ u(0)
\end{pmatrix} = M_E(p)^2 \begin{pmatrix} u(1) \\ u(0)
\end{pmatrix}
$$
and similarly
$$
\begin{pmatrix} u(p+1) \\ u(p)
\end{pmatrix} = M_E(p)^2 \begin{pmatrix} u(-p+1) \\ u(-p)
\end{pmatrix}.
$$
Moreover, \eqref{f.cayley-hamilton} implies
$$
M_E(p)^2 - \mathrm{Tr} M_E(p) \cdot M_E(p) + I = 0.
$$
Consequently, we have
\begin{equation}\label{f.gordoncase1}
\begin{pmatrix} u(2p+1) \\ u(2p)
\end{pmatrix} -  \mathrm{Tr} M_E(p) \begin{pmatrix} u(p+1) \\ u(p)
\end{pmatrix} + \begin{pmatrix} u(1) \\ u(0)
\end{pmatrix} = \begin{pmatrix} 0 \\ 0
\end{pmatrix}
\end{equation}
and
\begin{equation}\label{f.gordoncase2}
\begin{pmatrix} u(p+1) \\ u(p)
\end{pmatrix} -  \mathrm{Tr} M_E(p) \begin{pmatrix} u(1) \\ u(0)
\end{pmatrix} + \begin{pmatrix} u(-p+1) \\ u(-p)
\end{pmatrix} = \begin{pmatrix} 0 \\ 0
\end{pmatrix}.
\end{equation}
The assertion \eqref{f.gordonest} follows from
\eqref{f.gordoncase1} when $|\mathrm{Tr} M_E(p)| \le 1$ and it
follows from \eqref{f.gordoncase2} when $|\mathrm{Tr} M_E(p)| >
1$.
\end{proof}

The estimate \eqref{f.gordonest} can of course be used to exclude
the existence of decaying solutions. Notice that the energy $E$ does
not enter the argument. In particular, if the potential $V$ has
the required local periodicity for infinitely many values of $p$,
we have the estimate \eqref{f.gordonest} for infinitely many
values of $p$. This in turn shows that no $E$ can be an
eigenvalue. It is clear that one can perturb about this situation
a little bit and still deduce useful estimates. In light of this,
the following definition is natural.

\begin{definition}
A bounded potential $V : \Z \to \R$ is called a Gordon potential
if there are positive integers $q_k \to \infty$ such that
\begin{equation}\label{f.gordondef2}
\forall \, C > 0 : \lim_{k \to \infty} \max_{1 \le n \le q_k}
|V(n) - V(n \pm q_k)| C^{q_k} = 0.
\end{equation}
\end{definition}

\begin{lemma}\label{l.gordon2}
Suppose $V$ is a Gordon potential. Then, the operator $H = \Delta
+ V$ has purely continuous spectrum. More precisely, for every $E
\in \R$ and every solution $u$ of $H u = E u$, we have
\begin{equation}\label{f.gordonpotsolest}
\limsup_{|n| \to \infty} \left\| \begin{pmatrix} u(n+1)
\\ u(n) \end{pmatrix} \right\| \ge \frac12 \left\| \begin{pmatrix} u(1) \\ u(0)
\end{pmatrix} \right\|.
\end{equation}
\end{lemma}

\begin{proof}
By assumption, there is a sequence $q_k \to \infty$ such that
\eqref{f.gordondef2} holds. Given $E \in \R$, we consider a
solution $u$ of $H u = E u$ and, for every $k$, a solution $u_k$
of
$$
u_k(n+1) + u_k(n-1) + V_k(n) u_k(n) = E u_k(n)
$$
with $u_k(1) = u(1)$ and $u_k(0) = u(0)$, where $V_k$ is the
$q_k$-periodic potential that coincides with $V$ on the interval
$1 \le n \le q_k$.

It follows from Lemma~\ref{l.gordon} that $u_k$ satisfies the
estimate
$$
\max \left\{ \left\| \begin{pmatrix} u_k(-q_k+1) \\ u_k(-q_k)
\end{pmatrix} \right\| , \left\| \begin{pmatrix} u_k(q_k+1) \\ u_k(q_k)
\end{pmatrix} \right\| , \left\| \begin{pmatrix} u_k(2q_k+1) \\ u_k(2q_k)
\end{pmatrix} \right\| \right\} \ge \frac12 \left\| \begin{pmatrix} u(1) \\ u(0)
\end{pmatrix} \right\|.
$$
Since $V$ is very close to $V_k$ on the relevant interval and $u$
and $u_k$ have the same initial conditions, we expect that they
are close throughout the relevant interval and hence $u$ obeys a
similar estimate.

Let us make this observation explicit. Denote the transfer
matrices associated with $V_k$ by $M_{k,E}(n)$. We have
\begin{align*}
\max_{-q_k \le n \le 2q_k} \left\| \begin{pmatrix} u(n+1) \\
u(n) \end{pmatrix} - \begin{pmatrix} u_k(n+1) \\
u_k(n) \end{pmatrix}\right\| & \le \max_{-q_k \le n \le 2q_k}
\left\| M_E(n) - M_{k,E}(n) \right\| \left\| \begin{pmatrix} u(1)
\\ u(0) \end{pmatrix} \right\| \\
& \le 2 q_k C^{q_k} \max_{-q_k \le n \le 2q_k} |V(n) - V_k(n)| \left\|
\begin{pmatrix} u(1)
\\ u(0) \end{pmatrix} \right\|,
\end{align*}
which goes to zero by \eqref{f.gordondef2}.
\end{proof}

\begin{proof}[Proof of Theorem~\ref{t.avsim}]
Let $\alpha$ be Liouville and denote by $\frac{p_k}{q_k}$ the associated rational approximants. Consider arbitrary $\omega \in \T$ and $C > 0$. Since
$|\cos x - \cos y| \le |x - y|$, we have
$$
\max_{1 \le n \le q_k} |2 \lambda \cos (2 \pi (\omega + n \alpha)) - 2 \lambda \cos (2 \pi (\omega + (n \pm q_k) \alpha))| C^{q_k}
\le 4 \lambda \pi \, \mathrm{dist}(q_k \alpha , \Z) \, C^{q_k}
$$
and it follows from the Liouville condition that the right-hand side goes to zero as $k \to \infty$. Thus, $V(n) = 2 \lambda \cos (2 \pi (\omega + n \alpha))$ is a Gordon potential and hence $H^{\lambda,\alpha}_\omega$ has no eigenvalues by Lemma~\ref{l.gordon2}.
\end{proof}

\subsection{Resonant Phases and the Jitomirskaya-Simon Method}

In the previous subsection we saw how local translation symmetries
of the potential induce certain local translation symmetries of
solutions. While the latter are quite weak, they suffice to
exclude decay at infinity and hence square-summability. In this
section we briefly discuss an analogous study based on reflection
symmetries.

\begin{theorem}[Jitomirskaya-Simon 1994]\label{t.jitosim}
Given $\lambda$ and $\alpha$, there is a dense $G_\delta$ set of
$\omega$'s for which
$H_\omega^{\lambda,\alpha}$ has purely continuous spectrum.
\end{theorem}

\begin{definition}
A bounded potential $V : \Z \to \R$ is called a Jitomirskaya-Simon potential
if there are
$$
B > 4 \log (3 + 2 \|V\|_\infty)
$$
and integers $m_k \to \infty$ such that for every $k$,
$$
\sup_{n \in \Z} |V(2 m_k - n) - V(n)| < e^{-B m_k}.
$$
\end{definition}

\begin{lemma}
If $V$ is a Jitomirskaya-Simon potential and $E \in \R$, then
$$
u(n+1) + u(n-1) + V(n) u(n) = E u(n)
$$
has no non-zero $\ell^2$ solutions $u$.
\end{lemma}

\begin{proof}[Proof of Theorem~\ref{t.jitosim}]
Since $\cos$ is an even function, there is a dense set of $\omega$'s,
namely $\{k\alpha\}$, that have a center around which the associated
potential is symmetric. Take balls of suitable radius to generate open
sets which cover $\T$. The $\limsup$ of these balls is then a dense
$G_\delta$ set.
\end{proof}


\section{Localization at Supercritical Coupling}

In this section we discuss the main steps in the proof of
Theorem~\ref{t.jitoloc}.

\subsection{Existence of Generalized Eigenfunctions} We say that
$E \in \R$ is a generalized eigenvalue of $H$ if $Hu = Eu$ has a
non-trivial solution $u_E$, called the corresponding generalized
eigenfunction, satisfying
\begin{equation}\label{f.geneigenf}
|u_E(n)| \le C(1 + |n|)^\delta
\end{equation}
for suitable finite constants $C$ and $\delta$, and every $n \in
\Z$.

\begin{theorem}
{\rm (a)} Every generalized eigenvalue of $H$ belongs to
$\sigma(H)$.
\\
{\rm (b)} Fix $\delta > \frac12$ and some spectral measure $\mu$. Then, for $\mu$-almost every $E
\in \R$, there exists a generalized eigenfunction satisfying
\eqref{f.geneigenf}.
\\
{\rm (c)} The spectrum of $H$ is given by the closure of the set
of generalized eigenvalues of $H$.
\end{theorem}

\subsection{Solutions and Green's Function} For
$$
[n_1,n_2] = \{ n \in \Z : n_1 \le n \le n_2 \},
$$
denote by
$H_{[n_1,n_2]}$ the restriction of $H$ to this interval, that is,
$$
H_{[n_1,n_2]} = P_{[n_1,n_2]} H P_{[n_1,n_2]}^*,
$$
where $P_{[n_1,n_2]} : \ell^2(\Z) \to \ell^2([n_1,n_2])$ is the canonical projection
and $P_{[n_1,n_2]}^* : \ell^2([n_1,n_2]) \to \ell^2(\Z)$ is the canonical embedding.

Moreover, for
$E \not\in \sigma(H_{[n_1,n_2]})$ and $n,m \in [n_1,n_2]$, let
$$
G_{[n_1,n_2]}(n,m;E) = \langle \delta_n , \left( H_{[n_1,n_2]} - E
\right)^{-1} \delta_m \rangle.
$$
Then, the following formula holds.

\begin{lemma}\label{l.eigengreenlocal}
Suppose $n \in [n_1,n_2] \subset \Z$ and $u$ is a solution of the
difference equation $Hu=Eu$. If $E \not\in
\sigma(H_{[n_1,n_2]})$ and $n,m \in [n_1,n_2]$, then
$$
u(n) = - G_{[n_1,n_2]}(n,n_1;E) u(n_1 - 1) - G_{[n_1,n_2]}(n,n_2;E)
u(n_2 + 1).
$$
\end{lemma}

\begin{proof}
Since $u$ is a solution, we have
\begin{align*}
0 & = P_{[n_1,n_2]} (H - E) u \\
& = P_{[n_1,n_2]} (H - E) P_{[n_1,n_2]}^* P_{[n_1,n_2]} u +
P_{[n_1,n_2]} (H - E) P_{\Z \setminus [n_1,n_2]}^* P_{\Z \setminus
[n_1,n_2]} u,
\end{align*}
which in turn implies
$$
\left( H_{[n_1,n_2]} - E \right) (P_{[n_1,n_2]} u) = - \left(
u(n_1 - 1) \delta_{n_1} + u(n_2 + 1) \delta_{n_2} \right).
$$
Thus, with the inner product of $\ell^2([n_1,n_2])$, we find for
$n \in [n_1,n_2]$,
\begin{align*}
u(n) & = \langle \delta_n , (P_{[n_1,n_2]} u) \rangle \\
& = \left\langle \delta_n , \left( H_{[n_1,n_2]} - E \right)^{-1}
\left( H_{[n_1,n_2]} - E \right) (P_{[n_1,n_2]} u) \right\rangle \\
& = - \left\langle \delta_n , \left( H_{[n_1,n_2]} - E
\right)^{-1} \left( u(n_1 - 1) \delta_{n_1} + u(n_2 + 1)
\delta_{n_2} \right) \right\rangle \\
& = - u(n_1 - 1) \left\langle \delta_n , \left( H_{[n_1,n_2]} - E
\right)^{-1} \delta_{n_1} \right\rangle - u(n_2 + 1) \left\langle
\delta_n , \left( H_{[n_1,n_2]} - E \right)^{-1} \delta_{n_2}
\right\rangle
\end{align*}
as claimed.
\end{proof}

\subsection{Green's Function and Determinants}

Let
$$
P_k(\omega,E) = \det \left( (H_\omega^{\lambda,\alpha} - E)_{[0,k-1]} \right ).
$$
By Cramer's Rule, we have
for $n_1$, $n_2 = n_1 + k - 1$, and $n \in [n_1,n_2]$,
\begin{align*}
\left| G_{[n_1,n_2]}(n_1,n;E) \right| & = \left|
\frac{P_{n_2-n}(\omega + (n+1)\alpha,E)}{P_k(\omega + n_1
\alpha,E)}
\right|, \\
\left| G_{[n_1,n_2]}(n,n_2;E) \right| & = \left|
\frac{P_{n-n_1}(\omega + n_1 \alpha,E)}{P_k(\omega + n_1
\alpha,E)} \right|.
\end{align*}

\subsection{Determinants and Lyapunov Exponents}

We have
\begin{equation}\label{f.tmpolyrel}
M_E(k,\omega) = \begin{pmatrix} P_k(\omega,E) & - P_{k-1}(\omega +
\alpha, E) \\
P_{k-1}(\omega,E) & -P_{k-2} (\omega + \alpha , E) \end{pmatrix}.
\end{equation}
as can be checked by considering the degree, leading coefficient and zeros of
the entries of the transfer matrix, regarded as polynomials in $E$.

\begin{lemma}\label{l.jitolem0}
For every $E \in \R$ and $\varepsilon > 0$, there exists
$k(E,\varepsilon)$ such that
$$
|P_k(\omega,E)| < \exp\left( (\gamma(E) + \varepsilon) k \right)
$$
for every $k > k(E,\varepsilon)$ and every $\omega \in \T$.
\end{lemma}

\begin{proof}
This is a consequence of \eqref{f.tmpolyrel} and
the subadditive ergodic theorem for uniquely ergodic transformations.
\end{proof}

\begin{definition}
Fix $E \in \R$ and $\gamma \in \R$. A point $n \in \Z$ is
called $(\gamma,k)$-regular if there exists an interval
$[n_1,n_2]$, containing $n$ such that
\begin{align*}
(\mathrm{i}) & \quad n_2 = n_1 + k - 1, \\
(\mathrm{ii}) & \quad n \in [n_1,n_2], \\
(\mathrm{iii}) & \quad |n-n_i| > \frac{k}{5}, \\
(\mathrm{iv}) & \quad \left| G_{[n_1,n_2]}(n,n_i;E) \right| < \exp
\left( -\gamma |n-n_i| \right).
\end{align*}
Otherwise, $n$ is called $(\gamma,k)$-singular.
\end{definition}

The following central lemma describes the repulsion of singular clusters:

\begin{lemma}\label{l.jitokeylem}
For every $n \in \Z$, $\varepsilon > 0$, $\tau < 2$, there exists
$k_1 = k_1(\omega,\alpha,n,\varepsilon,\tau,E)$ such that for
every
$$
k \in \mathcal{K} = \left\{ k \in \Z_+ : \exists \, \tilde \omega \in \T \text{
with } |P_k(\tilde \omega,E)| \ge \frac{1}{\sqrt{2}} \, e^{k \gamma(E)}
\right\}
$$
with $k > k_1$, we have that
$$
m, n \text{ are both $(\gamma(E) - \varepsilon , k)$-singular and
} |m-n| > \frac{k+1}{2} \; \Rightarrow \; |m-n| > k^\tau.
$$
\end{lemma}

\begin{proof}[Proof of Theorem~\ref{t.jitoloc}.]
Let $E(\omega)$ be a generalized eigenvalue of
$H_\omega^{\lambda,\alpha}$ and denote the corresponding
generalized eigenfunction by $u_E$.
Notice that every point $n \in \Z$ with $u_E(n) \not=0$ is $(\gamma,k)$-singular for
$k > k_2 = k_2(E,\gamma,\omega,n)$.

Assume without loss of generality $u_E(0) \not= 0$ (otherwise replace zero by one).
Thus, by Lemma~\ref{l.jitokeylem}, if
$$
|n| > \max \left\{ k_1(\omega, \alpha, 0 , \varepsilon, 1.5 , E) , k_2(E,\gamma(E) -
\varepsilon , \omega , 0) \right\} + 1,
$$
the point $n$ is $(\gamma(E) - \varepsilon , k)$-regular for some
$k \in \{ |n|-1, |n|, |n|+1 \} \cap \mathcal{K} \not= \emptyset$,
since $0$ is $(\gamma(E) - \varepsilon , k)$-singular. Thus, there
exists an interval $[n_1 , n_2]$ of length $k$ containing $n$ such
that
$$
\frac{1}{5} (|n|-1) \le |n-n_i| \le \frac{4}{5} (|n|+1)
$$
and
$$
\left| G_{[n_1,n_2]}(n,n_i;E) \right| < e^{-(\gamma(E) -
\varepsilon) |n-n_i|}.
$$
From this and Lemma~\ref{l.eigengreenlocal}, we therefore see that
$$
|u_E(n)| \le 2 C(u_E) (2|n| + 1) e^{- \left( \frac{\gamma(E) -
\varepsilon}{5} \right) (|n|-1)}.
$$
By the uniform lower bound $\gamma(E) \ge \log \lambda$, this
implies exponential decay if $\varepsilon$ is chosen small enough.
\end{proof}

\begin{proof}[Sketch of the proof of Lemma~\ref{l.jitokeylem}.]
Assume that $m_1$ and $m_2$ are both $(\gamma(E) - \varepsilon ,
k)$-singular with
$$
d = m_2-m_1 > \frac{k+1}{2}.
$$
We set $n_i = m_i - \left\lfloor \frac{3}{4} k \right\rfloor$, $i =
1,2$.

It may be shown that there is a polynomial $Q_k$ of degree
$k$ such that
\begin{equation}\label{f.pkrep}
P_k(\tilde \omega) = Q_k\left( \cos \left( 2\pi (\tilde \omega + \frac{k-1}{2}
\alpha \right) \right).
\end{equation}
Let
$$
\omega_j = \begin{cases} \omega + \left( n_1 + \frac{k-1}{2} + j
\right) \alpha  , & j = 0,1,\ldots, \left\lfloor \frac{k+1}{2}
\right\rfloor - 1, \\ \omega + \left( n_2 + \frac{k-1}{2} + j -
\left\lfloor \frac{k+1}{2} \right\rfloor \right) \alpha  , & j =
\left\lfloor \frac{k+1}{2} \right\rfloor , \left\lfloor
\frac{k+1}{2} \right\rfloor + 1,\ldots, k. \end{cases}
$$
Lagrange interpolation then shows
\begin{equation}\label{f.jitoprooflagrange}
\left| Q_k(z) \right| = \left| \sum_{j=0}^k Q_k(\cos(2\pi\omega_j)
\frac{\prod_{l\not= j} (z - \cos (2\pi \omega_l)}{\prod_{l\not= j}
(\cos (2\pi \omega_j) - \cos (2\pi \omega_l))} \right|.
\end{equation}
By $(\gamma(E) - \varepsilon ,
k)$-singularity and Lemma~\ref{l.jitolem0}, we have for $k$ sufficiently large,
\begin{equation}\label{f.jitoproofclaim1}
|Q_k(\cos(2\pi\omega_j) | < \exp \left( \frac{k}{8}(\gamma(E) -
\varepsilon) \right), \quad j = 0,1,\ldots,k.
\end{equation}

Using that the Diophantine and non-resonance assumptions, one may show that if $d < k^\tau$ for some $\tau < 2$,
we have for large $k$,
\begin{equation}\label{f.jitoproofclaim2}
\frac{|\prod_{l\not= j} (z - \cos (2\pi \omega_l)|}{|\prod_{l\not=
j} (\cos (2\pi \omega_j) - \cos (2\pi \omega_l))|} \le
\exp\left(\frac{k \varepsilon}{16} \right) \; \text{ for $z \in
[-1,1]$, } 0 \le j \le k.
\end{equation}

Given $\tau < 2$, consider $k \in \mathcal{K}$ large enough and $\tilde
\omega$ with
$$
|P_k(\tilde \omega)| \ge \frac{1}{\sqrt{2}} \, e^{k \gamma(E)}.
$$
But assuming $d < k^\tau$, we also have the following upper bound,
$$
|P_k(\tilde \omega)| \le (k+1) \exp \left( \frac{k}{8}(\gamma(E) -
\varepsilon) \right) \exp\left(\frac{k \varepsilon}{16} \right),
$$
which follows from \eqref{f.pkrep}--\eqref{f.jitoproofclaim2} for
$z = \cos \left( 2\pi (\tilde \omega + \frac{k-1}{2} \alpha
\right)$. This contradiction shows that $d < k^\tau$ is
impossible.
\end{proof}


\section{Cantor Spectrum via Aubry Duality}
In this section we present Puig's proof of the striking fact that
localization for the operator family
$\{H^{\lambda,\alpha}_\omega\}_{\omega \in \T}$, as established in
the previous lecture, implies via Aubry duality and reducibility
Cantor spectrum for the dual family
$\{H^{\lambda^{-1},\alpha}_\omega\}_{\omega \in \T}$. Once one has
shown Cantor spectrum for $0 < \lambda < 1$, Aubry duality applied
again yields Cantor spectrum for $\lambda > 1$. That is, we will
show how Theorem~\ref{t.puigcantor} follows from Jitomirskaya's
localization result.

Consider the equations
\begin{equation}\label{f.amoeve}
u(n+1) + u(n-1) + 2 \lambda \cos(2 \pi n \alpha) u(n) = E u(n),
\end{equation}
\begin{equation}\label{f.damoeve}
u(n+1) + u(n-1) + 2 \lambda^{-1} \cos(2 \pi (\omega + n \alpha)) u(n) =
(\lambda^{-1} E) u(n)
\end{equation}
and recall from Section~2 how exponentially decaying
eigenfunctions and highly regular quasi-periodic solutions are
related to each other via Aubry duality:

\begin{lemma}\label{l.puiglem2}
{\rm (a)} Suppose $u$ is an exponentially decaying solution of
\eqref{f.amoeve}. Consider its Fourier series
$$
\hat u(\omega) = \sum_{m \in \Z} u(m) e^{2 \pi i m \omega}.
$$
Then, $\hat u$ is real-analytic on $\T$, it extends analytically
to a strip, and the sequence $\tilde u(n) = \hat u(\omega + n \alpha)$ is a
solution of \eqref{f.damoeve}.

{\rm (b)} Conversely, suppose $u$ is a solution of
\eqref{f.damoeve} with $\omega = 0$ of the form $u(n) = g(n\alpha)$ for some
real-analytic function $g$ on $\T$. Consider the Fourier series
$$
g(\omega) = \sum_{n \in \Z} \hat g (n) e^{2 \pi i n \omega}.
$$
Then, the sequence $\{\hat g (n)\}$ is an exponentially decaying
solution of \eqref{f.amoeve}.
\end{lemma}

\begin{proof}
(a) Since $u$ is exponentially decaying, $\hat u$ extends to a
function that is analytic in a neighborhood of the unit circle
$\T$. The other statement was shown in Section~2.

(b) Since $u$ is a solution of \eqref{f.damoeve} with $\omega = 0$ and we have $u(m)
= g(m\alpha)$, we have
$$
g((m+1) \alpha) + g((m-1)\alpha) + 2 \lambda^{-1} \cos(2 \pi m
\alpha) g(m\alpha) = (\lambda^{-1} E) g(m\alpha).
$$
Rewriting this in terms of the Fourier expansion, we find
\begin{align*}
\sum_{n \in \Z} \hat g (n) e^{2 \pi i n (m+1) \alpha} + \sum_{n
\in \Z} \hat g (n) e^{2 \pi i n (m-1) \alpha} + \lambda^{-1} &
\left( e^{- 2 \pi i m \alpha} + e^{2 \pi i m \alpha} \right)
\sum_{n \in \Z} \hat g (n) e^{2 \pi i n m \alpha} \\
& = (\lambda^{-1} E) \sum_{n \in \Z} \hat g (n) e^{2 \pi i n m
\alpha}.
\end{align*}
It follows that
\begin{align*}
& \sum_{n \in \Z} E \hat g (n) e^{2 \pi i n m \alpha}\\ & = \sum_{n
\in \Z} \left( \lambda \hat g (n) e^{2 \pi i n (m+1) \alpha} +
\lambda \hat g (n) e^{2 \pi i n (m-1) \alpha} + \left( e^{- 2 \pi
i m \alpha} + e^{2 \pi i m \alpha} \right)  \hat g (n) e^{2 \pi i
n m \alpha} \right) \\
& = \sum_{n \in \Z} \left( 2 \lambda \cos(2 \pi n \alpha) \hat g
(n) + \hat g (n-1) + \hat g (n+1) \right) e^{2 \pi i n m \alpha}
\end{align*}
Since $g$ is real-analytic on $\T$ with analytic extension to a
strip, it follows that the Fourier coefficients of $g$ decay
exponentially and satisfy the difference equation
\eqref{f.amoeve}.
\end{proof}

Next we use the information provided by the previous
lemma to reduce the situation at hand to constant coefficients. We prove a general
statement to this effect:

\begin{lemma}\label{l.puiglem3}
Let $\alpha \in \T$ be Diophantine and suppose $A : \T \to
\mathrm{SL}(2,\R)$ is a real-analytic map, with analytic extension
to the strip $|\Im \omega| < \delta$ for some $\delta > 0$. Assume
that there is a non-vanishing real-analytic map $v : \T \to \R^2$
with analytic extension to the same strip $|\Im \omega| < \delta$
such that
$$
v(\omega + \alpha) = A(\omega) v(\omega) \quad \text{ for every }
\omega \in \T.
$$
Then, there are a real number $c$ and a real-analytic map $B : \T
\to  \mathrm{SL}(2,\R)$ with analytic extension to the strip $|\Im
\omega| < \delta$ such that with
\begin{equation}\label{f.puigcdef}
C = \begin{pmatrix} 1 & c \\ 0 & 1 \end{pmatrix},
\end{equation}
we have
\begin{equation}\label{f.redtoconst}
B(\omega + \alpha)^{-1} A(\omega) B(\omega) = C \quad \text{ for
every } \omega \in \T.
\end{equation}
\end{lemma}

\begin{proof}
Since $v$ does not vanish, $d(\omega) = v_1(\omega)^2 +
v_2(\omega)^2$ is strictly positive and hence we can define
$$
B_1(\omega) = \begin{pmatrix} v_1(\omega) & -
\frac{v_2(\omega)}{d(\omega)} \\ v_2(\omega) &
\frac{v_1(\omega)}{d(\omega)} \end{pmatrix} \in \mathrm{SL}(2,\R)
$$
for $\omega \in \T$. We have
\begin{equation}\label{f.puig1}
A(\omega) B_1 (\omega) = \begin{pmatrix} v_1(\omega + \alpha) &
\ast \\ v_2(\omega + \alpha) & \ast
\end{pmatrix} \in \mathrm{SL}(2,\R)
\end{equation}
and hence
$$
A(\omega) B_1(\omega) = B_1 (\omega + \alpha) \tilde C(\omega)
$$
with
$$
\tilde C(\omega) =
\begin{pmatrix} 1 & \tilde c(\omega) \\ 0 & 1
\end{pmatrix},
$$
where $\tilde c : \T \to \R$ is analytic. Indeed, by
\eqref{f.puig1} the first column of $\tilde C(\omega)$ is
determined and then its $(2,2)$ entry must be one since $\tilde
C(\omega) = B_1 (\omega + \alpha)^{-1} A(\omega) B_1(\omega) \in
\mathrm{SL}(2,\R)$. Now let
$$
c = \int_\T \tilde c(\omega) \, d\omega.
$$
and define the matrix $C$ as in \eqref{f.puigcdef}.

We claim that we can find $b : \T \to \R$ analytic (with analytic
extension to a strip) such that
\begin{equation}\label{f.cohomeq}
b(\omega + \alpha) - b(\omega) = \tilde c(\omega) - c \quad \text{
for every } \omega \in \T.
\end{equation}
Indeed, expand both sides of the hypothetical identity
\eqref{f.cohomeq} in Fourier series:
$$
\sum_{k \in \Z} b_k e^{2 \pi i (\omega + \alpha) k} - \sum_{k \in
\Z} b_k e^{2 \pi i \omega k} = \sum_{k \in \Z} \tilde c_k e^{2 \pi
i \omega k} - c.
$$
Since we have $\tilde c_0 = c$, the $k=0$ terms disappear on both
sides and hence all we need to do is to require
$$
b_k (e^{2 \pi i \alpha k} - 1) = \tilde c_k \quad \text{ for every
} k \in \Z \setminus \{ 0 \}.
$$
In other words, if we set $b_0 = 0$ and
$$
b_k  = \frac{\tilde c_k}{e^{2 \pi i \alpha k} - 1} \quad \text{
for every } k \in \Z \setminus \{ 0 \},
$$
then
$$
b(\omega) = \sum_{k \in \Z} b_k e^{2 \pi i \omega k}
$$
satisfies \eqref{f.cohomeq}. Since $\tilde c (\cdot)$ has an
analytic extension to a strip, the coefficients $\tilde c_k$ decay
exponentially. On the other hand, the Diophantine condition which
$\alpha$ satisfies ensures that the coefficients $b_k$ decay
exponentially as well and hence $b(\cdot)$ is real-analytic with
an extension to the same open strip.

Setting
$$
B_2(\omega) = \begin{pmatrix} 1 & b(\omega)
\\ 0 & 1 \end{pmatrix} \in
\mathrm{SL}(2,\R),
$$
and using \eqref{f.cohomeq}, we find
\begin{align*}
B_2(\omega + \alpha)^{-1} \tilde C(\omega) B_2 (\omega) & =
\begin{pmatrix} 1 & - b(\omega + \alpha)
\\ 0 & 1 \end{pmatrix}
\begin{pmatrix} 1 & \tilde c(\omega) \\ 0 & 1
\end{pmatrix}
\begin{pmatrix} 1 & b(\omega)
\\ 0 & 1 \end{pmatrix} \\
& = \begin{pmatrix} 1 & - b(\omega + \alpha)
\\ 0 & 1 \end{pmatrix}
\begin{pmatrix} 1 & b(\omega) + \tilde c(\omega) \\ 0 & 1
\end{pmatrix} \\
& = \begin{pmatrix} 1 & b(\omega) + \tilde c(\omega) - b(\omega +
\alpha)
\\ 0 & 1 \end{pmatrix} \\
& = \begin{pmatrix} 1 & c
\\ 0 & 1 \end{pmatrix} \\
& = C
\end{align*}
for every $\omega \in \T$. Thus, setting $B(\omega) = B_1 (\omega)
B_2(\omega)$, we obtain \eqref{f.redtoconst}.
\end{proof}

\begin{proof}[Proof of Theorem~\ref{t.puigcantor}.]
Consider first a coupling constant $\lambda > 1$. We have seen
above that Aubry duality maps the energy $E$ to the dual energy
$\lambda^{-1} E$. We will establish below that if $E$ is an
eigenvalue of $H^{\lambda,\alpha}_0$, then the dual energy
$\lambda^{-1} E$ is an endpoint of a gap of the spectrum of
$H^{\lambda^{-1},\alpha}_0$. Since we already know that for
Diophantine $\alpha$, $H^{\lambda,\alpha}_0$ has pure point
spectrum, we can consider energies belonging to the countable
dense set of eigenvalues. It then follows that the dual energies
are all endpoints of gaps and hence the gaps are dense because the
spectra are just related by uniform scaling.

Let us implement this strategy. Consider an eigenvalue $E$ of
$H^{\lambda,\alpha}_0$ and a corresponding exponentially decaying
eigenfunction. Then, Lemma~\ref{l.puiglem2} yields the
real-analytic function $\hat u$, which has an analytic extension
to a strip, and a quasi-periodic solution of the dual difference
equation at the dual energy. Using this as input to
Lemma~\ref{l.puiglem3}, we then obtain that
$$
A(\omega) = \begin{pmatrix} \lambda^{-1} E - 2 \lambda^{-1} \cos(2
\pi \omega) & - 1 \\ 1 & 0 \end{pmatrix}
$$
may be analytically conjugated via $B(\cdot)$ to the constant
$$
C = \begin{pmatrix} 1 & c \\ 0 & 1 \end{pmatrix}.
$$

Let us show that $c \not= 0$. Assume to the contrary $c = 0$. Then,
$A(\omega) = B(\omega + \alpha) B(\omega)^{-1}$
for every $\omega \in \T$ and therefore, all solutions of
\eqref{f.damoeve} are analytically quasi-periodic! Indeed,
\begin{align*}
\begin{pmatrix} u(n) \\ u(n-1) \end{pmatrix} & = A(\omega + (n-1)
\alpha) \begin{pmatrix} u(n-1) \\ u(n-2) \end{pmatrix} \\
& = \cdots \\
& = A(\omega + (n-1) \alpha) \times \cdots \times A(\omega) \begin{pmatrix} u(0) \\
u(-1) \end{pmatrix} \\
& = B(\omega + n \alpha) B(\omega)^{-1} \begin{pmatrix} u(0) \\
u(-1) \end{pmatrix},
\end{align*}
that is,
$$
u(n) = \left\langle \begin{pmatrix} 1 \\ 0 \end{pmatrix} , B(\omega + n \alpha) B(\omega)^{-1} \begin{pmatrix} u(0) \\
u(-1) \end{pmatrix} \right\rangle,
$$
and hence $u(n) = g(n\alpha)$ with a real-analytic function $g$ on
$\T$. Now consider two linearly independent solutions of
\eqref{f.damoeve} and associate with them via
Lemma~\ref{l.puiglem2} the corresponding exponentially decaying
solutions of the dual equation \eqref{f.amoeve}. They must be
linearly independent too, which yields the desired contradiction
since by constancy of the Wronskian there cannot be two linearly
independent exponentially decaying solutions. This contradiction shows
$c \not= 0$.

Let us now perturb the energy and consider
$$
\tilde A(\omega) = \begin{pmatrix} (\lambda^{-1} E + \lambda^{-1}
\delta) - 2 \lambda^{-1} \cos(2 \pi \omega) & - 1 \\ 1 & 0
\end{pmatrix} = A(\omega) + \begin{pmatrix} \lambda^{-1}
\delta & 0 \\ 0 & 0 \end{pmatrix}.
$$
One can show that there is $\delta_0 > 0$ such that
\begin{equation}\label{f.puigncg}
0 < | \delta | < \delta_0 \text{ and } \delta c < 0 \quad
\Rightarrow \quad \lambda^{-1} E + \lambda^{-1} \delta \not\in
\sigma(H^{\lambda^{-1},\alpha}_0).
\end{equation}

Since the $E$'s in question are dense in
$\sigma(H^{\lambda,\alpha}_0)$, Aubry duality shows that the
$\lambda^{-1} E$'s in question are dense in
$\sigma(H^{\lambda^{-1},\alpha}_0)$. By \eqref{f.puigncg} all
these energies are endpoints of gaps of
$\sigma(H^{\lambda^{-1},\alpha}_0)$. Thus,
$\sigma(H^{\lambda^{-1},\alpha}_0)$ does not contain an interval.
Recall that by general principles,
$\sigma(H^{\lambda^{-1},\alpha}_0)$ is closed and does not contain
isolated points. Consequently, $\Sigma^{\lambda^{-1},\alpha} =
\sigma(H^{\lambda^{-1},\alpha}_0)$
is a Cantor set. Then, by Aubry duality again,
$\Sigma^{\lambda,\alpha} = \lambda \Sigma^{\lambda^{-1},\alpha}$ is a
Cantor set, too. Putting everything together,
it follows that for every Diophantine $\alpha$ and every $\lambda \in
(0,\infty) \setminus \{ 1 \}$, $\Sigma^{\lambda,\alpha}$ is a Cantor set.
\end{proof}


\section{Liouville Frequencies}

So far we have shown the main theorems, that is, an identification
of the almost sure spectral type and the Cantor structure of the
spectrum, for Diophantine frequencies $\alpha$ and non-critical
coupling constants $\lambda$. In this final section we discuss the
case of Liouville frequencies in more depth. They have already
been considered when we showed that there are never any
eigenvalues and hence the expected localization at super-critical
coupling fails in these cases.

As a consequence of this non-result, we cannot deduce using Aubry
duality that the spectrum at sub-critical coupling is almost surely
purely absolutely continuous nor that it is a Cantor set. Thus, a
different approach is required to prove these statements, which turn
out to be indeed true. Since Liouville numbers are very well
approximated by rational numbers, it is natural to try and prove the
expected statements by periodic approximation. This is also the
philosophy that is implemented with the help of Gordon's lemma to
prove the absence of eigenvalues. Periodic operators always have
purely absolutely continuous spectrum, so one needs a method to push
this through to the quasi-periodic limit. Moreover, the spectra of
periodic operators have many gaps and a suitable continuity statement
for spectra could feasibly identify a dense set of gaps in the
spectrum of the quasi-periodic operator if the rate of approximation
is sufficiently good.

We have the following pair of theorems. For simplicity, we work
with the notion of a Liouville number as introduced earlier.
Stronger results are known but the proofs are (much) more
difficult.

\begin{theorem}[Choi-Elliott-Yui 1990]
Suppose $\alpha$ is Liouville. Then, for every $\lambda > 0$,
$\Sigma^{\lambda,\alpha}$ is a Cantor set.
\end{theorem}

\begin{theorem}[Avila-Damanik 2008]
Suppose $\alpha$ is Liouville. Then, for every $0 < \lambda < 1$ and
almost every $\omega \in \T$, $H^{\lambda,\alpha}_\omega$ has purely
absolutely continuous spectrum.
\end{theorem}

We will only make these results plausible by describing the main
ideas and tools used in the proofs. As pointed out above, it will
be essential to obtain a good understanding of the periodic
approximants to $H^{\lambda,\alpha}_\omega$, which arise when
$\alpha$ is replaced by a rational number close to it. So, from
now on, we will drop the requirement that $\alpha$ is irrational
and instead consider $H^{\lambda,\alpha}_\omega$ for $\lambda > 0$
and $\alpha,\omega \in \T$. Given such $\lambda$ and $\alpha$, we
define
$$
\Sigma^{\lambda,\alpha} = \bigcup_{\omega \in \T}
\sigma(H^{\lambda,\alpha}_\omega).
$$
For irrational $\alpha$, this definition coincides with the one we had
previously, since the spectrum of $H^{\lambda,\alpha}_\omega$ is
$\omega$-independent in this case.

We first discuss Cantor spectrum. By Aubry duality, it suffices to
consider $0 < \lambda \le 1$. This will sometimes be relevant
below. The first question we address is the continuity of
$\Sigma^{\lambda,\alpha}$ as a function of $\alpha$.

\begin{lemma}[Avron-van Mouche-Simon 1990]
For every $\lambda > 0$, there exists $\delta > 0$ such that if
$|\alpha - \alpha'| < \delta$, then
$$
\mathrm{dist}_H \left(\Sigma^{\lambda,\alpha} ,
  \Sigma^{\lambda,\alpha'}\right) \le 6 \left( 2 \lambda |\alpha -
  \alpha'| \right)^{1/2}.
$$
\end{lemma}

\begin{figure}
\includegraphics[width=0.9\textwidth]{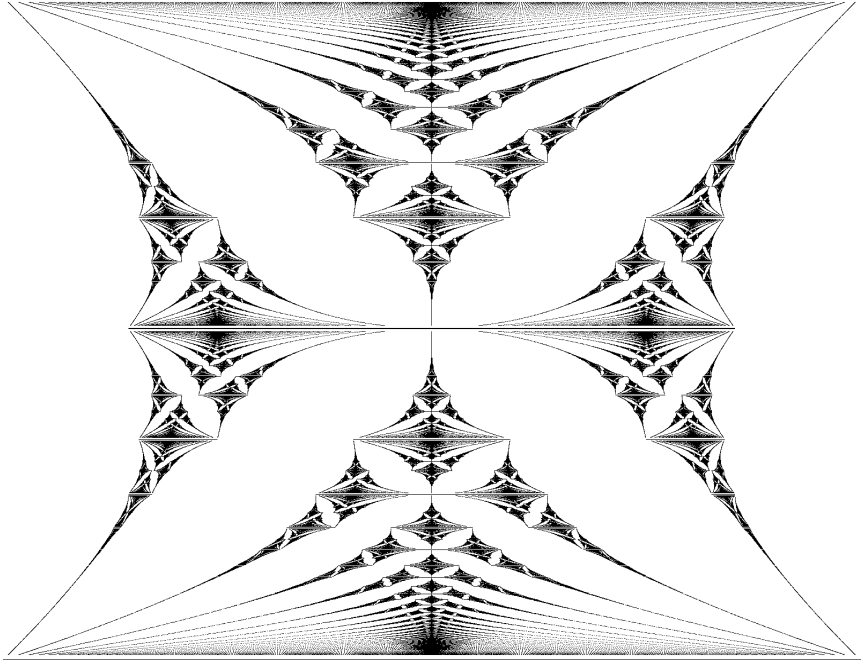} \par
\caption{The Hofstadter butterfly: The $x$-axis corresponds to the
  energy $E$ and the $y$-axis corresponds to the frequency $\alpha$.}
\label{fig:hof}
\end{figure}

Here, $\mathrm{dist}_H (\Sigma^{\lambda,\alpha} ,
\Sigma^{\lambda,\alpha'})$ denotes the Hausdorff distance between
$\Sigma^{\lambda,\alpha}$ and $\Sigma^{\lambda,\alpha'}$. The proof
uses test functions with a linear cut-off. This continuity result
shows that any gap of $\Sigma^{\lambda,p/q}$ corresponds to a gap of
$\Sigma^{\lambda,\alpha}$ if its length is larger than $6 \left( 2
  \lambda |\alpha - \alpha'| \right)^{1/2}$. Figure~\ref{fig:hof}
shows the sets $\Sigma^{1,p/q}$ for $0 \le p/q \le 1$ with $q \le
50$. This plot, which has a beautiful self-similar structure, is known
as the \textit{Hofstadter butterfly}.

The following result gives quite detailed information about the number
and size of the gaps of $\Sigma^{\lambda,p/q}$. The result should be
compared with the plot of the Hofstadter butterfly.

\begin{lemma}[Choi-Elliott-Yui 1990]
Suppose $0 < \lambda \le 1$ and $p,q$ are coprime positive
integers. Write $q = 2m+1$ or $q = 2m+2$. Then, $\Sigma^{\lambda,p/q}$
has exactly $2m$ gaps. Each of these gaps is of length at least
$\lambda^m 8^{-q}$.
\end{lemma}

This produces many gaps in $\Sigma^{\lambda,\alpha}$ when $\alpha$ is
Liouville. Indeed, choose a sequence $p_k/q_k \to \alpha$ with
$p_k,q_k$ coprime and $|\alpha - p_k/q_k| < k^{-q_k}$. Then, for $k$
large enough, there are at least $q_k - 2$ many gaps in
$\Sigma^{\lambda,\alpha}$. It remains to show that these gaps are
dense. To this end, the following observation does the job.

\begin{lemma}[Choi-Elliott-Yui 1990]
Suppose $0 < \lambda \le 1$ and $p,q$ are coprime positive
integers. For every gap of $\Sigma^{\lambda,p/q}$, there is another
gap of $\Sigma^{\lambda,p/q}$ within distance $8\pi / q$.
\end{lemma}

Let us now discuss how to prove purely absolutely continuous spectrum
at sub-critical coupling for Liouville frequencies. The spectral type
is much less robust under perturbations than, for example, the
spectrum. Simple approximation of $H^{\lambda,\alpha}_\omega$ by
$H^{\lambda,p/q}_\omega$ does not look promising at first. It turns
out that robustness improves if one averages over $\omega$!

Recall that the spectral measure $d\mu_\omega$ associated with
$H^{\lambda,\alpha}_\omega$ and the state $\delta_0 \in \ell^2(\Z)$ is
the Borel probability measure with the property that
$$
\langle \delta_0 , g(H^{\lambda,\alpha}_\omega) \delta_0 \rangle =
\int_\R g(E) \, d\mu_\omega(E)
$$
for all bounded, Borel measurable functions $g$. The density of states
measure $dN$ is formally given by
$$
dN(E) = \int_\T d\mu_\omega(E) \, d\omega.
$$
More precisely, it is the Borel probability measure that satisfies
$$
\int_\T \langle \delta_0 , g(H^{\lambda,\alpha}_\omega) \delta_0
\rangle \, d\omega = \int_\R g(E) \, dN(E)
$$
for all bounded, Borel measurable functions $g$. Consider the Lebesgue
decomposition of these measures. Denote the Radon-Nikodym derivative
of the absolutely continuous component of $d\mu_\omega$ (resp., $dN$)
by $\mu_{\omega , \mathrm{ac}}(E)$ (resp., $N_\mathrm{ac}(E)$).

\begin{lemma}[Kotani 1997]
For almost every $E \in \{ E : \gamma(E) = 0 \}$, we have
$$
N_\mathrm{ac}(E) = \int_\T \mu_{\omega , \mathrm{ac}}(E) \, d\omega.
$$
\end{lemma}

\begin{theorem}[Kotani 1997]
Suppose that the Lyapunov exponent vanishes on $\Sigma^{\lambda,\alpha}$ and
$$
\int_\R N_\mathrm{ac}(E) \, dE = 1.
$$
Then, $H^{\lambda,\alpha}_\omega$ has purely absolutely continuous
spectrum for almost every $\omega \in \T$.
\end{theorem}

Since the periodic approximants have purely absolutely continuous
spectrum, the associated density of states measures are absolutely
continuous so we want to take a limit of these quantities as $p_k/q_k
\to \alpha$. What about the Lyapunov exponents? It is well known that
the Lyapunov exponent vanishes on the spectrum in the periodic
case. For the almost Mathieu operator, we have the following result:

\begin{theorem}[Bourgain-Jitomirskaya 2002]
If $\lambda > 0$ and $\alpha$ is irrational, then
$\gamma(E) = \max \{ \log \lambda , 0 \}$ for every $E \in
\Sigma^{\lambda,\alpha}$.
\end{theorem}

This shows that for $0 < \lambda \le 1$, the Lyapunov exponent
vanishes on the spectrum and we can attempt to apply Kotani's theorem
in this coupling regime. We mention in passing that for $\lambda = 1$,
the spectrum has zero Lebesgue measure and hence $dN$ is purely
singular. For $0 < \lambda < 1$, on the other hand, and a Liouville
$\alpha$ it is indeed possible to show that the irrational case is so
well approximated by the rational case that
$\int_\R N_\mathrm{ac}(E) \, dE = 1$. Kotani's theorem then yields the result.

\bigskip

\noindent\textbf{Acknowledgments.} The author would like to thank
the members of the research group of Michael Baake at the
University of Bielefeld for the kind hospitality and the CRC~701
for generous financial support.

\end{document}